# Is the Halting problem effectively solvable non-algorithmically, and is the Gödel sentence in NP, but not in P?

Bhupinder Singh Anand[1]

We consider the thesis that an arithmetical relation, which holds for any, given, assignment of natural numbers to its free variables, is Turing-decidable if, and only if, it is the standard representation of a PA-provable formula. We show that, classically, such a thesis is, both, unverifiable and irrefutable, and, that it implies the Turing Thesis is false; that Gödel's arithmetical predicate $R(x)$, treated as a Boolean function, is in the complexity class NP, but not in P; and that the Halting problem is effectively solvable, albeit not algorithmically.

## 1. Introduction

Classical theory holds that:

(*a*) Every Turing-computable function $F$ is partial recursive[2], and, if $F$ is total[3], then $F$ is recursive ([Me64], p233, Corollary 5.13).

---

[1] The author is an independent scholar. E-mail: re@alixcomsi.com; anandb@vsnl.com. Postal address: 32, Agarwal House, D Road, Churchgate, Mumbai - 400 020, INDIA. Tel: +91 (22) 2281 3353. Fax: +91 (22) 2209 5091.

[2] Classically ([Me64], p120-121, p214), a partial function $F$ of $n$ arguments is called partial recursive if, and only if, $F$ can be obtained from the initial functions (zero function), projection functions, and successor function (of classical recursive function theory) by means of substitution, recursion and the classical, unrestricted, $\mu$-operator. $F$ is said to come from $G$ by means of the unrestricted $\mu$-operator, where $G(x_1, ..., x_n, y)$ is recursive, if, and only if, $F(x_1, ..., x_n) = \mu y(G(x_1, ..., x_n, y) = 0)$, where $\mu y(G(x_1, ..., x_n, y) = 0)$ is the least number $k$ (if such exists) such that, if $0 =< i =< k$, $G(x_1, ..., x_n, i)$ exists and is not 0, and $G(x_1, ..., x_n, k) = 0$. We note that, classically, $F$ may not be defined for certain $n$-tuples; in particular, for those $n$-tuples $(x_1, ..., x_n)$ for which there is no $y$ such that $G(x_1, ..., x_n, y) = 0$.

(*b*) Every partial recursive function is Turing-computable ([Me64], p237, Corollary 5.15).

From this, it concludes that the following, essentially unverifiable[4] but refutable[5], theses are equivalent ([Me64], p237:

*Church's Thesis*: A number-theoretic function is effectively computable if, and only if, it is recursive ([Me64], p227).

*Turing's Thesis*: A number-theoretic function is effectively computable if, and only if, it is Turing-computable ([BBJ03], p33).

In this paper, we now show how, under another intuitively unobjectionable thesis, Turing's Thesis does not hold, and consider some consequences.

We note, first, that, even classically, the above equivalence, informally referred to as CT, does not hold strictly, and needs further qualification. The following argument highlights this, where *F* is any number-theoretic function:

(*i*) Assume Church's Thesis. Then:

If *F* is Turing-computable, then, by (*a*), it is partial recursive. If *F* is total, then it is both recursive ([Me64], p227) and, by our assumption, effectively computable.

---

[3] We define a number-theoretic function, or relation, as total if, and only if, it is effectively computable, or effectively decidable, respectively, for any, given, set of natural number values assigned to its free variables. We define a number-theoretic function, or relation, as partial otherwise. We define a partial number theoretic function, or relation, as effectively computable, or decidable, respectively, if, and only if, it is effectively computable, or decidable, respectively, for any, given, set of values, assigned to its free variables, for which it is defined.

[4] The two theses are essentially unverifiable in classical theory since the notion of 'effective computability' is intuitive, and not defined formally.

[5] Demonstration of a number-theoretic function that is effectively computable, but not recursive, would falsify Church's Thesis; similarly, demonstration of a number-theoretic function that is effectively computable, but not Turing-computable, would falsify Turing's Thesis.

If *F* is effectively computable, then, by our assumption, it is recursive. Hence, by definition, it is partial recursive and, by (*b*), Turing-computable.

(*ii*) Assume Turing's Thesis. Then:

If *F* is recursive, it is partial recursive and, by (*b*), Turing-computable. Hence, by our assumption, *F* is effectively computable.

If *F* is effectively computable, then, by our assumption, it is Turing-computable. Hence, by (*a*), it is partial recursive and, if *F* is total, then it is recursive.

Thus, for CT to be a strict equivalence, every partial recursive function would need to be effectively extendable as a total partial recursive function.

Now, it follows from Turing's reasoning [Tu36] that such an assumption is, in fact, inconsistent with classical theory.

Thus, in his seminal paper on computable numbers [Tu36], Turing considers the Halting problem, which can be expressed as the query:

*Halting problem for T*: Given a Turing machine T, can one effectively decide, given any instantaneous description *alpha*, whether or not there is a computation of T beginning with *alpha* ([Me64], p256)?

Turing then shows that the Halting problem is unsolvable by a Turing machine.

Since a function is Turing-computable if, and only if, it is partially Markov-computable ([Me64], p233, Corollary 5.13 & p237, Corollary 5.15), it is essentially unverifiable algorithmically whether, or not, a Turing machine that computes a random, $n$-ary, number-theoretic function will halt classically, without going into a non-terminating loop, on every $n$-ary sequence of natural numbers, for which it is defined, as input, where:

*Non-terminating loop*: A non-terminating loop is defined as any repetition of the instantaneous tape description[6] of a Turing machine during a computation.

It follows that we cannot extend every partial recursive function as a total recursive function by defining its undefined values algorithmically.

In this paper we now show that, under an intuitively unobjectionable thesis, every partial recursive function can, nevertheless, be effectively extended as a total function. Since such an extended definition is not algorithmic, the extended function is not partial recursive. However, since it is effectively computable, the Turing Thesis is false under the assumed thesis.

## 2. The Provability Thesis

Now, in Theorem VI of his seminal paper [Go31a] on undecidable propositions, Gödel constructs a primitive recursive relation that can be shown to be true for all assignments of natural number values to its free variables. Hence, treated as a Boolean function, it is a total function that is Turing-computable[7]. However, none of its representations in first order Peano Arithmetic, PA, are PA-provable.

More precisely, Gödel constructs - in an intuitionistically unobjectionable manner - an arithmetical predicate, say $[R(x)]$[8], that is unprovable in PA, but which is Tarskian-true

---

[6] "An instantaneous tape description describes the condition of the machine and the tape at a given moment. When read from left to right, the tape symbols in the description represent the symbols on the tape at the moment. The internal state $q_s$ in the description is the internal state of the machine at the moment, and the tape symbol occurring immediately to the right of $q_s$ in the tape description represents the symbol being scanned by the machine at the moment." ([Me64], p230, footnote 1).

[7] Since every total recursive function is Turing-computable ([Me64], p237, Corollary 5.15).

[8] We use square brackets to indicate that the expression inside them refers to a particular syntactical string of symbols that is effectively verifiable as a well-defined formula of a specified formal system.

under the standard interpretation of PA since, given any natural number $n$, [$R(n)$] is PA-provable.[9]

Further, in his Theorem VII ([Go31a], p29), Gödel shows that every primitive recursive relation is arithmetical. The result can be extended to show that every recursive relation, treated as a Boolean function, is representable in PA, and is, therefore, also arithmetical ([Me64], p131-134, Propositions 3.23 & 3.24).

Taken together, the above suggest the possibility that:

> *Arithmetical Uncomputability Thesis*: There is a total arithmetical function that is not Turing-computable.

Unlike the Church and Turing Theses, this thesis is, classically, irrefutable. This follows since, by Turing's Halting argument, Turing-computability is essentially unverifiable.

However, we note that:

> **Corollary (*a*)**: The *Arithmetical Uncomputability Thesis* implies that the Turing Thesis is false.[10]

---

[9] This follows from Gödel's argument that, if $r$ is the Gödel-number of the formula [$R(x)$] in his Arithmetic P, then the P-formula, [$R(n)$], whose Gödel-number is $Sb(r, 17|Z(n))$, is provable for any given natural number $n$ ([Go31a], p26).

[10] However, this does not imply that the Church Thesis, too, is false. Under the assumptions of this essay, even though a number-theoretic function is Turing-computable if, and only if, it is partial recursive, the equivalence between Church's Thesis and Turing's Thesis may not hold.

Thus, whereas Gödel's primitive recursive relation $\sim xB(Sb(p\ 19|Z(p)))$ - where $p$ is the Gödel-number of the formula [$R(x, y)$] that represents $\sim xB(Sb(y\ 19|Z(y)))$ (cf. [Go31a], p24, def. 8.1) in a Peano Arithmetic - would be Turing decidable, it's instantiationally equivalent arithmetical relation $R(x, p)$ would not be Turing decidable (thus, in a narrow sense, falsifying Turing's Thesis)!

In other words, the truth of Gödel's unprovable-but-true arithmetical proposition, [($Ax$)$R(x, p)$], under the standard interpretation, would follow verifiably from the Turing-decidability of the primitive recursive relation $\sim xB(Sb(p\ 19|Z(p)))$, but not from that of $R(x, p)$; the latter, like the Halting predicate, would be Turing-undecidable under the *Provability* thesis (and, so, any Turing machine that computes its

**Corollary (*b*)**: The *Arithmetical Uncomputability Thesis* implies that Gödel's arithmetical predicate $R(x)$, treated as a Boolean function, is in the complexity class NP[11], but not in P.[12]

In what follows, we replace the Arithmetical Uncomputability Thesis by another, intuitively unobjectionable, Provability Thesis, and use this to construct a total arithmetical relation that, treated as a Boolean function, is not Turing-computable.

We consider, first, the following thesis:

*Decidability Thesis*: An arithmetical relation, which holds for any, given, assignments of natural numbers to its free variables, is Turing-decidable if, and only if, it is the standard representation of a PA-provable formula.

This too, is, classically, both unverifiable and irrefutable, since, again, Turing's Halting argument implies that we cannot assume that there is always an algorithm that will verify whether an arithmetical relation, which holds for all assignments of natural numbers to its free variables, is Turing-decidable.

However, the Decidability Thesis is of particular interest, since it echoes some implicitly held beliefs in interpretations of computational theory. For instance, in an arXived paper ([CCS01], v2), Calude et al hold that:

---

characteristic function would not halt on every natural number $x$, but would loop non-terminatingly for some natural number $n$).

[11] Informally, a problem is assigned to the complexity class P (deterministic polynomial-time) if the number of steps needed to solve it is bounded by some power of the problem's size. A problem is assigned to the complexity class NP (nondeterministic polynomial-time) if it permits a nondeterministic solution and the number of steps needed to verify the solution is bounded by some power of the problem's size. (Cf. [We05])

[12] This follows since, if $R(x)$, treated as a Boolean function, were in P, then it would be Turing-computable.

Classically, there are two equivalent ways to look at the mathematical notion of proof: logical, as a finite sequence of sentences strictly obeying some axioms and inference rules, and computational, as a specific type of computation. Indeed, from a proof given as a sequence of sentences one can easily construct a Turing machine producing that sequence as the result of some finite computation and, conversely, given a machine computing a proof we can just print all sentences produced during the computation and arrange them into a sequence.

In other words, the authors seem to hold - in the sense of the Decidability Thesis - that Turing-computability of a 'proof', in the case of an arithmetical proposition, is equivalent to provability of its representation in PA.

Now, the Decidability Thesis can be expressed more precisely as the assertion that:

*Provability Thesis*: When computing an arithmetical relation $F(x_1, ..., x_n)$, treated as a Boolean function[13], a classical[14] Turing machine T halts, and returns a value that interprets as 'true', without going into a non-terminating loop, on every $n$-ary sequence of natural numbers as input if, and only if, $[F(x_1, ..., x_n)]$ is a PA-provable formula.

---

[13] In other words, the Turing machine computes the characteristic function, say $C(x_1, ..., x_n)$, of $F(x_1, ..., x_n)$, which is defined as follows ([Me64], p119):

$C(x_1, ..., x_n) = 0$ if $F(x_1, ..., x_n)$ is true;
$C(x_1, ..., x_n) = 1$ if $F(x_1, ..., x_n)$ is false.

[14] We take Mendelson [Me64], Boolos et al [BBJ03], and Rogers [Ro87], as representative, in the areas that they cover, of standard expositions of classical, first order, logic and of effective computability (in particular, of standard Peano Arithmetic and of classical Turing-computability).

Further, we note that:

> *Effective Looping oracle*: Any Turing machine can be provided with an auxiliary infinite tape ([Ro87], p130) to effectively recognise a non-terminating looping situation; it simply records[15] every instantaneous tape description at the execution of each machine instruction on the auxiliary tape, and compares the current instantaneous tape description with the record. If an instantaneous tape description is repeated, it can be meta-programmed to abort the impending non-terminating loop, and return a meta-symbol indicating self-termination.

We now show that the Provability Thesis implies that we can effectively determine, albeit non-algorithmically, whether, or not, T will halt on a given input. In other words, the Provability Thesis implies that the Halting problem for Turing machines is effectively solvable - but not by a Turing machine.

## 3. An effective, non-algorithmic, solution of the Halting problem

We detail this argument below.

> **Meta-theorem 1**: The Provability Thesis implies that every partial recursive number-theoretic function $F(x_1, ..., x_n)$ can be effectively extended as a total function.
>
> **Proof**: We assume that $F$ is obtained from the recursive function $G$ by means of the unrestricted $\mu$-operator, so that $F(x_1, ..., x_n) = \mu y(G(x_1, ..., x_n, y) = 0)$.
>
> If $[H(x_1, ..., x_n, y)]$ expresses $\sim(G(x_1, ..., x_n, y) = 0)$ in PA, we consider the PA-provability, and Tarskian-truth (cf. [Me64], p49), in the standard interpretation M of PA, of the formula $[H(a_1, ..., a_n, y)]$ for a given sequence of numerals $<[a_1], ..., [a_n]>$ of PA.

---

[15] It is convenient to visualise the tape of such a Turing machine as that of a two-dimensional virtual-teleprinter, which maintains a copy of every instantaneous tape description in a random-access memory during a computation.

Thus:

(*a*) Let $Q_1$ be the meta-assertion that $[H(a_1, ..., a_n, y)]$ is not classically true in M. Since $G(a_1, ..., a_n, y)$ is recursive, it follows that there is some finite $k$ such that any Turing machine $T_1(y)$ that computes $G(a_1, ..., a_n, y)$ will halt and return the value 0 for $y = k$.

(*b*) Next, let $Q_2$ be the meta-assertion that $[H(a_1, ..., a_n, y)]$ is classically true in M, but that there is no Turing machine that computes the arithmetical function $H(a_1, ..., a_n, y)$ as a Boolean function without going into a non-terminating loop.

Since $G(a_1, ..., a_n, y)$ is recursive, it follows that there is some finite $k'$ such that any Turing machine $T_2(y)$ that computes the arithmetical function $H(a_1, ..., a_n, y)$ as a Boolean function will halt, as its auxiliary tape will return the symbol for self-termination at the occurrence of a non-terminating loop for $y = k'$.

(*c*) Finally, let $Q_3$ be the meta-assertion that $[H(a_1, ..., a_n, y)]$ is classically true in M, and that any Turing machine that computes the arithmetical function $H(a_1, ..., a_n, y)$ as a Boolean function does not enter into any non-terminating loop.

Now, if we assume the Provability Thesis, then it follows that $[H(a_1, ..., a_n, y)]$ is PA-provable. Let $h$ be the Gödel-number of $[H(a_1, ..., a_n, y)]$. We consider, then, Gödel's primitive recursive number-theoretic relation $xBy$ ([Go31a], p22, def. 45), which holds in M if, and only if, $x$ is the Gödel-number of a proof sequence in PA for the PA-formula whose Gödel-number is $y$. It follows that there is some finite $k''$ such that any Turing machine $T_3(y)$, which computes the characteristic function[16] of $xBh$, will halt and return the value 0 for $x = k''$[17].

Since $Q_1$, $Q_2$, and $Q_3$ are mutually exclusive and exhaustive[18], it follows that, when run simultaneously[19] over the sequence 1, 2, 3, ... of values for $y$, one of the parallel trio $\{T_1(y) // T_2(y) // T_3(y)\}$ of Turing machines will always halt for some finite value of $y$. If $T_1(y)$ halts at $y = k$, and returns the value 0, we define $F'(a_1, ..., a_n) = F(a_1, ..., a_n)$. If $T_2(y)$ halts and returns the symbol for self-termination, or if $T_3(y)$ halts, we define $F'(a_1, ..., a_n) = 0$.

Hence, under the given hypotheses, there is always an effective extension of every partial recursive function, $F(x_1, ..., x_n)$, as a total function, $F'(x_1, ..., x_n)$.¶

---

[16] "If $R(x_1, ..., x_n)$ is a relation, then the characteristic function $C_n(x_1, ..., x_n)$ is defined as follows:

$C_n(x_1, ..., x_n) = 0$ if $R(x_1, ..., x_n)$ is true, and $C_n(x_1, ..., x_n) = 1$ if $R(x_1, ..., x_n)$ is false." ([Me64], p119).

[17] We assume that such a machine can be effectively meta-programmed to recognise, and proceed to the next instantaneous tape description whenever it encounters, a non-terminating loop.

[18] They correspond to the instances where a classical Turing machine that computes the arithmetical relation $H(a_1, ..., a_n, y)$, treated as a Boolean function, will halt for some $y$, go into a non-terminating loop for some $y$, or not halt for any $y$, respectively.

[19] This concept is essentially that of parallel computing, where the action of one machine can influence the action of another unpredictably, without human intervention. Since classical Turing machines are necessarily sequential, such a procedure cannot be defined as a classical Turing machine. Hodges remarks [HA00] that the possibility of parallel machines being essentially different from his Logical Computing Machines does not (arguably) appear to have been considered by Turing:

"... Another source may lie in Turing's definition of an 'oracle-machine' which is a Turing machine allowed at certain points to 'consult the oracle'. Such a machine is not purely mechanical: it is like the 'choice-machine' defined in (Turing 1936-7) which at certain points allows for human choices to be made. Turing used the word 'machine' for entities which are only partially mechanical in operation, reserving the term 'automatic machine' for those which are purely mechanical. Copeland appears to imagine that when Turing describes the oracle-machine definition as giving a 'new type of machine', he is defining a new type of automatic machine. On the contrary, Turing is defining something only partially mechanical.

To take this point further, it is worth noting that the expression 'purely mechanical process' enters into Turing's definitive statement of the Church-Turing thesis, which comes as an opening section to (Turing 1939), and that Turing goes on: 'understanding by a purely mechanical process one which could be carried out by a machine'. In the subsequent discussion the word 'machine' is used to mean 'Turing machine'. There is no evidence that Turing had any concept of a purely mechanical 'machine' of any kind other than encapsulated by the Turing machine definition."

Since the Halting problem is not Turing-solvable classically, it follows that:

**Corollary 1.1**: The Provability Thesis implies that the parallel trio $\{T_1(y) \mathbin{/\mkern-5mu/} T_2(y) \mathbin{/\mkern-5mu/} T_3(y)\}$ of Turing machines is not a Turing machine.

**Corollary 1.2**: The Provability Thesis implies that the classical Halting problem for Turing machines is effectively solvable, but not algorithmically.

Further, since the Provability Thesis implies the Arithmetical Unprovability Thesis, it follows that:

**Corollary 1.3**: The Provability Thesis implies that the Turing Thesis is false.

**Corollary 1.4**: The Provability Thesis implies that Gödel's arithmetical predicate $R(x)$, treated as a Boolean function, is in the complexity class NP, but not in P.

(*Updated: Wednesday 8$^{nd}$ June 2005 4:12:46 AM IST by re@alixcomsi.com*)